\documentclass[a4paper]{article}
\usepackage{amsmath}
\usepackage{amsfonts}
\usepackage{amssymb}         
\usepackage{amsopn}
\usepackage{theorem}          
\usepackage{latexsym}
\usepackage{graphicx}        
\usepackage{mathrsfs}		
\usepackage{psfrag}
\usepackage[noend]{algorithmic}
\usepackage{algorithm}
\usepackage{color}
\usepackage{verbatim}
\usepackage{url}
\usepackage{multirow}
\usepackage{setspace}
\usepackage{booktabs,longtable,rotating}
\usepackage{afterpage}

\theoremstyle{changebreak}                

 {{\nopagebreak\hspace*{\fill}$\Box$\par\vspace{12pt}}}

\newcommand{\tS}{@{$\;$}}

\begin{document}
\thispagestyle{empty}
\begin{center} 

{\LARGE A Local Branching Heuristic for MINLPs}
\par \bigskip
{\sc Giacomo Nannicini${}^1$, Pietro Belotti${}^2$, Leo Liberti${}^1$}
\par \bigskip
\begin{minipage}{15cm}
\begin{flushleft}
{\small
\begin{itemize}
\item[${}^1$] {\it LIX, \'Ecole Polytechnique, F-91128 Palaiseau,
France} \\ Email:\url{{giacomon,liberti}@lix.polytechnique.fr}
\item[${}^2$] {\it Lehigh University, 200 West Packer Avenue,
  Bethlehem, PA 18015} \\ Email:\url{belotti@lehigh.edu}
\end{itemize}
}
\end{flushleft}
\end{minipage}
\par \medskip \today
\end{center}
\par \bigskip

\begin{abstract}
Local branching is an improvement heuristic, developed within the
context of branch-and-bound algorithms for MILPs, which has proved to
be very effective in practice. For the binary case, it is based on
defining a neighbourhood of the current incumbent solution by allowing
only a few binary variables to flip their value, through the addition
of a local branching constraint. The neighbourhood is then explored
with a branch-and-bound solver. We propose a local branching scheme
for (nonconvex) MINLPs which is based on iteratively solving MILPs and
NLPs. Preliminary computational experiments show that this approach is
able to improve the incumbent solution on the majority of the test
instances, requiring only a short CPU time. Moreover, we provide
algorithmic ideas for a primal heuristic whose purpose is to find a
first feasible solution, based on the same scheme.
\end{abstract}

\section{Introduction}
\label{s:intro}
Local branching was introduced by Fischetti and Lodi
\cite{localbranching} as a primal heuristic for Mixed-Integer Linear
Programs (MILPs) within the context of a Branch-and-Bound (BB)
algorithm. It is sometimes referred to as an {\it improvement}
heuristic, in that it aims at improving the primal bound by finding a
better feasible solution, starting from the incumbent (i.e.~the best
known feasible solution). The natural setting for local branching is
within problems with binary variables, although it has also been
extended to the case of general integer variables
\cite{localbranching}. The idea is as follows: whenever a new
incumbent is found, a new problem is solved, which has the same
feasible region and objective of the original problem, but with the
addition of a local branching constraint, whose purpose is to allow
only a given number of binary variables to change their value with
respect to the incumbent. This new problem, which we call the {\it
  local branching problem}, defines a neighbourhood which is explored
by employing a BB algorithm. The computational time required to solve
the local branching problem is typically very small, as the local
branching constraint greatly reduces the feasible region. Therefore,
it is some kind of local search. Moreover, the local branching problem
need not be solved to optimality: when employed as heuristic for
difficult problems, it is typically solved with a maximum running
time, or until a better solution is found. Computational experiments
have shown that this simple idea is often able to improve the
incumbent on a large number of real-world test problems, thus being
practically useful to reduce the running time of a BB algorithm by
providing a better bound. Clearly, this paradigm can be directly
applied to BB methods for nonconvex Mixed-Integer Nonlinear Programs
(MINLPs): the only difference is that the local branching problem
should be solved by employing a BB algorithm for MINLPs. However,
branch-and-bound methods for this class of problems are in practice
significantly slower than in the linear case, because branching can
occur on integer and continuous variables
\cite{tawarmalani1,belottibbt}, a convexification refinement step is
applied, and large continuous Nonlinear Programs (NLPs) are solved at
some nodes. Therefore, the exploration of the neighbourhood may
require more time, losing effectiveness.

In this paper, we propose a local branching scheme for nonconvex
MINLPs, which is based on repeatedly solving a sequence of
limited-size MILPs and NLPs. Each of these problems has a different
purpose. We use a NLP to estimate the descent direction of the
original objective function. Then, we solve a MILP on the
convexification of the feasible region with a local branching
constraint, to enforce integral feasibility. Finally, we fix the
integer variables and try to satisfy the original constraint via a
NLP. In case of failure, we cut off the computed solution from the
MILP, and iterate the algorithm. As each solved problem has a smaller
feasible region or involves fewer variables with respect to the
original MINLP, this approach is fast. We report preliminary
computational experiments to assess its usefulness. We also propose
ideas for a primal heuristic whose purpose is to obtain a first
feasible solution at the root node of the BB tree. Future research
will include the implementation of these heuristics within an existing
solver for nonconvex MINLPs, to test them in practice on a larger
number of instances, and a full version of this preliminary paper.

\section{Theoretical Background}
\label{s:background}
Consider the following mathematical program:
\begin{equation}
  \left. \begin{array}{rr\tS c\tS l} 
    \min & f(x) && \\
    \forall j \in M & g_j(x) &\le& 0\\
    \forall i \in N & x_i^L \le x_i &\le& x_i^U \\
    \forall j \in N_I & x_j &\in& \mathbb{Z}, 
  \end{array} \right\} \mathcal{P}
  \label{standard}
\end{equation}
where $f$ and $g$ are possibly nonconvex functions, $n = |N|$ is the
number of variables, and $x = (x_i)_{i \in N}$ is the vector of
variables. This general type of problem has applications in several
fields \cite{biegler,floudas,lln2}. Difficulties arise from the integrality
of some of the variables, as well as nonconvexities. Solution methods
typically require that the functions $f$ and $g$ are factorable, that
is, they can be expressed as $\sum_i \prod_j h_{ij}(x)$
\cite{tawarmalani1}. When the function $f$ is linear and the $g_j$'s
are affine, $\mathcal{P}$ is a MILP, for which efficient
Branch-and-bound or Branch-and-Cut methods have been developed
\cite{nemhauser,wolsey}. Commercial codes (e.g.~\cite{cplex110}) are
often able to solve MILPs with thousands of variables in reasonable
time. Branch-and-Bound methods for MINLPs attempt to closely mimick
their MILP counterparts, but many difficulties have to be overcome. In
particular, obtaining lower bounds is not straightforward. The
continuous relaxation of each subproblem may be a nonconvex NLP, whose
global optimum is difficult to find. One possibility is to compute a
convexification of the feasible region of the problem, so that lower
bounds can be easily computed. In the following, we will assume that
the Branch-and-Bound method that we address computes a linear
convexification of the original problem; that is, the objective
function $f$ and all constraints $g_j$ are replaced by suitable linear
terms which underestimate and overestimate the original functions over
all the feasible region. The accuracy of the convexification greatly
depends on the variable bounds. If the interval over which a variable
is defined is large, the convexification of functions which contain
that variable may be a very loose estimation of the original
functions. For this reason, branching can also occur on continuous
variables, so as to reduce the variable bounds. This improves the
quality of the convexification. Moreover, bounds can be tightened by
applying various techniques, such as Feasibility Based Bound
Tightening (FBBT), Optimality Based Bound Tightening (OBBT), etc.~(see
\cite{belottibbt}). We note that a good incumbent solution not only
provides a better upper bound, but also allows for the propagation of
tighter bounds through expression tree based bound tightening
techniques. Moreover, OBBT benefits from a better upper
bound. Therefore, finding good feasible solutions is doubly important
for BB algorithms for MINLPs.

Local Branching \cite{localbranching} is an improvement heuristic for
BB algorithms which relies on exploring a neighbourhood of the
incumbent, looking for a better solution. For problems with only
continuous and binary variables, the neighbourhood is defined by
adding a {\it local branching constraint} to the original problem,
obtaining the local branching problem. Let $B \subset N_I$ be the set
of binary variables, $0 < k \in \mathbb{N}$, and let $\bar{x}$ be any
feasible solution; then the local branching constraint is:
\begin{equation}
\label{eq:lb}
\sum_{i \in B: \bar{x}_i = 1} (1 - x_i) + \sum_{i \in B: \bar{x}_i =
  0} x_i \le k.
\end{equation}
This constraint has the effect of allowing only $k$ binary variables
to flip their value from $0$ to $1$ or vice versa. Typically, $k$ is a
small value; experiments in \cite{localbranching} suggest $k \approx
10$. As a consequence, the number of feasible solutions of the local
branching problem is very small, and an efficient BB code requires
little time to find its optimal solution. The heuristic was proposed
and applied as a primal heuristic for BB algorithms for MILPs; it has
also been used in conjunction with other metaheuristics, such as
Variable Neighbourhood Search (VNS) \cite{vnsfoundations}, both in the
context of MILPs \cite{vnslb}, and of nonconvex MINLPs
\cite{recipe}. In particular, the latter paper reports very good
results over a large collection of possibly nonconvex MINLPs by
applying an iterative exploration of the solution space defining
neighbourhoods of increasing size (in the spirit of VNS), where the
neighbourhood for binary variables is defined through
(\ref{eq:lb}). As the majority of the test instances have binary
variables, this turns out to be effective, and further motivates our
interest for local branching in the context of nonconvex
MINLPs. However, in \cite{recipe} the local branching neighbourhood is
explored by means of a solver for convex MINLPs \cite{bonmin,minlpbb},
i.e.~MINLPs whose continuous relaxation is convex. In this case, the
solvers are employed as a heuristic. For this purpose, \cite{bonmin}
suggests using a BB algorithm choosing the branching variables via NLP
strong branching, which implies solving several NLPs at each node of
the BB tree. As a result, the solution of the local branching problem
may be slow. Within the context of a BB algorithm for nonconvex
MINLPs, a local branching heuristic should be as fast as possible.

\section{Local Branching for MINLPs}
\label{s:lb}
Branch-and-Bound solvers for nonconvex MINLPs are slower than for
MILPs. There are several reasons for this. First, the convexification
of the problem may be computationally expensive. The convexification
is carried out at the root node, but it is typically refined at
various stages of the optimization process. This is also true for the
bound tightening phase. Second, branching can occur on integer and
continuous variables, therefore there is an overhead because more
possible branching variables have to be dealt with. Third, continuous
NLPs are solved at some nodes of the BB tree. Moreover, available
software for nonconvex MINLPs does not have the same reliability and
speed as solvers for MILPs, which have been tested and improved for
almost 20 years. All these difficulties motivate our idea for a local
branching scheme which does not employ a solver for nonconvex MINLPs
to solve the local branching problem.

Let $\bar{x}$ be the incumbent which we want to improve. Let
$\bar{\cal P}$ be the linear relaxation of ${\cal P}$
with the addition of the local branching constraint (\ref{eq:lb}). Let
${\cal Q}$ be the continuous relaxation of ${\cal P}$, i.e., ${\cal
  P}$ with no integrality constraints, and let $\bar{\cal Q}$ be
${\cal Q}$ with the additional constraint (\ref{eq:lb}). A naive
approach would be to solve $\bar{\mathcal{P}}$ using a MILP solver,
and then, if the solution obtained is not feasible with respect to the
original constraints of $\mathcal{P}$, employ a local NLP solver
fixing the integer variables to regain feasibility. Two problems
arise. First, the convexification of the objective function of ${\cal
  P}$ may be very different from the original objective
function. Hence, optimizing with respect to the convexified objective
function could deteriorate the objective value. Second, the solution
of $\bar{\cal P}$ is likely to be an infeasible point with respect to
the original constraints of ${\cal P}$. We would like to find a point
which is as close as possible to the feasible region, so that
constraint feasibility can be regained by modifying the continuous
variables only, and keeping the integer variables fixed. To do so, we
solve the continuous relaxation $\bar{\cal Q}$ using a local NLP
solver. This yields a point $x'$ such that $f(x') \le f(\bar{x})$,
since we have relaxed integrality. Moreover, $x'$ is feasible with
respect to the constraints of ${\cal P}$, although it is typically not
integral feasible. We use $x'$ to estimate the descent direction of the
original objective function, i.e.~to indentify the region in which a
better incumbent could be found. Let $\bar{P}$ be the feasible region
of $\bar{\cal P}$. We find an integral feasible point by employing a
MILP solver on the problem:
\begin{equation}
  \label{eq:milp}
  \min_{x \in \bar{P}}\|x-x'\|_{\ell}.
\end{equation}
If $\ell = 1$ or $\ell = \infty$, (\ref{eq:milp}) is a MILP. If $\ell
= 2$, it can be solved as a Mixed-Integer Quadratic Program (MIQP). In
the following, we assume $\ell = 1$. Let $x'' = \arg \min_{x \in
  \bar{P}}\|x-x'\|_{\ell}$. By solving (\ref{eq:milp}), we hopefully find an
integral feasible point which is near $x'$, hence it is likely that
$x''$ is almost feasible and improves the objective value. In the
following step we fix the integer variables of $x''$, and solve ${\cal
  P}$ with a local NLP solver with starting point $x''$. We obtain a
new point $x^\ast$. If $x^\ast$ is feasible for ${\cal P}$ and
$f(x^\ast) < f(\bar{x})$, we have a new incumbent, and the algorithm
terminates with success. Otherwise, we append the constraint
$LB_{rev}(x^\ast)$:
\begin{equation}
  \label{eq:reverselb}
  \sum_{i \in B: x^\ast_i = 1} (1 - x_i) + \sum_{i \in B: x^\ast_i = 0}
  x_i \ge 1
\end{equation}
to $\bar{P}$, and iterate the algorithm. (\ref{eq:reverselb})
avoids finding a solution with the same values on the binary variables
as $x^\ast$. This way, at each iteration we find different
solutions. 

\begin{algorithm}[htb]
\begin{algorithmic}
\STATE Initialization: $stop \leftarrow {\tt false}$
\STATE Solve $\bar{\cal Q}$ with a local NLP solver, obtaining point
$x'$
\WHILE{$\neg stop$}
\STATE Solve $\min_{x \in \bar{P}}\|x-x'\|_{1}$ with a MILP solver,
obtaining point $x''$
\STATE Solve ${\cal P}$ with a local NLP solver and initial point
$x''$, keeping the integer variables fixed, obtaining point $x^\ast$
\IF{($x^\ast$ is not feasible for ${\cal P}$) $\lor $
  ($f(\bar{x}) \le f(x^\ast)$)}
\STATE Append $LB_{rev}(x^\ast)$ to $\bar{P}$
\ELSE
\STATE $stop \leftarrow {\tt true}$
\ENDIF
\ENDWHILE
\STATE {\bf return} $x^\ast$
\end{algorithmic}
\caption{Local Branching Heuristic for MINLPs}
\label{alg:lb}
\end{algorithm}

We give a description of our algorithm in
Algorithm~\ref{alg:lb}. Although Algorithm~\ref{alg:lb} may iterate
until an improved incumbent is found or one of employed solvers fails,
additional stopping criteria can be used, such as a maximum CPU time
or a maximum number of iterations. Trivially, if $|B|$ is the number
of binary variables, and $k$ is the rhs of (\ref{eq:lb}),
Algorithm~\ref{alg:lb} will stop after at most $\binom{|B|}{k}$
iterations, returning either an improved incumbent, or no
solution. This follows from the fact that, after each iteration, one
realization of the vector of binary variables is excluded from the set
of feasible solutions to (\ref{eq:milp}) (through the addition of
$LB_{rev}(x^\ast)$), and there are at most $\binom{|B|}{k}$ possible
combinations.



Algorithm~\ref{alg:lb} employs a MILP solver and a local NLP solver
only, therefore it does not rely on BB nonconvex MINLP solvers, which
would typically slow down the local branching heuristic. However, a
nonconvex MINLP solver guarantees to find an improved incumbent within
the neighbourhood defined by the local branching constraint, if one
exists. In this case, we are trading reliability for speed. In the
context of a BB software for nonconvex MINLPs, heuristics are supposed
to be fast, therefore this approach finds application.

\section{Computational Experiments}
In this section we provide preliminary computational experiments. The
heuristic was implemented with the AMPL scripting language in order to
test if it is able to find improved solutions, so as to simulate its
behaviour when integrated within a MINLP solver. We used {\tt couenne}
\cite{couenne} to obtain the convexification of the problems. As MILP
solver, we employed Cplex 11.0 \cite{cplex110}, whereas the local NLP
solver is {\tt ipopt} \cite{ipopt}. The tests were run on one core of
an Intel Centrino Duo clocked at 1.06 Ghz, on a machine with 1.5 GB
RAM. The right hand side $k$ of the local branching constraint was
computed as
\begin{equation*}
k = \min(15,\max(1,|B|/2)),
\end{equation*}
where $B$ is the set of binary variables. Algorithm~\ref{alg:lb} was
terminated after 10 iterations of the main loop. Cplex was run with
default parameters and maximum running time of 2 seconds, whereas {\tt
  ipopt} was run with the options {\tt expect\_infeasible\_problem,
  start\_with\_resto}. To test the heuristic, we ran {\tt couenne} for
10 minutes on instances with both binary and continuous variables
taken from the MINLPLib
(\url{http://www.gamsworld.com/minlp/minlplib.htm}), and recorded the
first feasible solutions found, up to a maximum of two. Then we
applied Algorithm~\ref{alg:lb} on each of them, trying to find a
better incumbent. Results are reported in Table~\ref{tab:results}. For
each instance, we record the objective value of each initial feasible
solution tested; then we report the iteration of
Algorithm~\ref{alg:lb} at which we found the first improved solution,
the value of the new incumbent, and the total required CPU time. We
also record the best solution found during the 10 iterations. Note
that, since we used an AMPL script, no data is shared between the
solvers, and each time the problem data structures have to be
initalized; a more clever implementation, integrated within the code
of a MINLP solver, is likely to be faster. Moreover, we did not tune
the parameters of the solvers. Therefore, we believe that running time
can be further reduced.

\begin{table}
\centering
\begin{tabular}{l|r|rrr|rrr|}
  & {\sc Initial} & \multicolumn{3}{c|}{\sc First Imp. Solution} &
  \multicolumn{3}{c|}{\sc Best Solution} \\
  {\sc Instance} & {\sc Solution} & {\sc It.} & {\sc Objective} & {\sc
    Time}
  & {\sc It.} & {\sc Objective} & {\sc Time}\\
  \hline
  {\tt csched1} & -28438.6 & 1 & -30639.3 & 0.316 & 1 & -30639.3 & 0.316\\
                & -29779.8 & 1 & -30639.3 & 0.292 & 1 & -30639.3 & 0.292\\
  {\tt csched1a}& -29719.5 & 1 & -30430.2 & 0.108 & 1 & -30430.2 & 0.108\\
                & -29862.4 & 1 & -30430.2 & 0.136 & 1 & -30430.2 & 0.136\\
  {\tt csched2} & -123261 & 1 & -135365 & 3.600 & 10 & -141523 & 23.777\\
                & -128347 & 1 & -137722 & 2.952 & 2 & -149076 & 5.296\\
  {\tt csched2a}& -139073 & 1 & -142403 & 2.560 & 7 & -143793 & 18.867\\
                & -151353 & 1 & -155252 & 3.908 & 4 & -155977 & 14.284\\
  {\tt elf}     & 4.03665 & 1 & 2.57999 & 0.156 & 10 & 2.47666 & 1.808\\
                & 2.04448 & 6 & 1.61999 & 0.164 & 10 & 1.40666 & 1.912\\
  {\tt eniplac} & -132010 & - & - & - & 1 & -131648 & 1.700 \\
                & -132117 & - & - & - & 1 & -131648 & 1.688 \\
  {\tt enpro48} & 189132 & 9 & 188887 & 20.393 & 9 & 188887 & 20.393 \\
  {\tt enpro56} & 276551 & 1 & 266094 & 2.480 & 1 & 266094 & 2.480 \\
                & 275296 & 1 & 266094 & 2.476 & 1 & 266094 & 2.476 \\
  {\tt ex1233}  & 161022 & - & - & - & 6 & 178588 & 5.220 \\
                & 155522 & - & - & - & 7 & 155522 & 60.351 \\
  {\tt fo7}     & 46.9636 & 1 & 46.2517 & 2.192 & 4 & 39.6749 & 9.112\\
                & 30.4694 & 2 & 28.9003 & 4.324 & 8 & 28.5415 & 17.161\\
  {\tt fo7\_2}  & 41.3952 & - & - & - & 2 & 41.3952 & 6.408 \\
                & 30.9608 & 1 & 28.0644 & 2.232 & 2 & 21.6761 & 4.356\\
  {\tt fo8}     & 39.2582 & - & - & - & 1 & 39.62582 & 2.216 \\
  {\tt fo9}     & 42.7099 & 1 & 42.0448 & 2.336 & 1 & 42.0448 & 2.336\\
                & 41.1547 & 1 & 39.5119 & 2.320 & 2 & 37.1786 & 4.536\\
  {\tt m7}      & 202.098 & 1 & 150.357 & 2.348 & 3 & 126.337 & 7.144\\
                & 175.142 & 1 & 144.505 & 2.244 & 3 & 136.905 & 6.568\\
  {\tt o7}      & 167.586 & 8 & 159.605 & 17.865 & 8 & 159.605 & 17.865\\
                & 161.379 & 6 & 158.269 & 13.508 & 6 & 158.269 & 13.508\\
  {\tt o7\_2}   & 161.337 & - & - & - & 5 & 161.337 & 11.244\\
                & 127.366 & 3 & 116.946 & 6.500 & 3 & 116.946 & 6.500\\
  {\tt ravem}   & 295020 & 1 & 269590 & 1.236 & 1 & 269590 & 1.236 \\
                & 283851 & 1 & 269590 & 1.256 & 1 & 269590 & 1.256 \\
  {\tt st\_e36} & 0 & 1 & -2 & 0.416 & 1 & -2 & 0.416 \\
                & -1 & 1 & -2 & 0.380 & 1 & -2 & 0.380 \\
  {\tt water4}  & 1645.76 & 1 & 1022.47 & 2.456 & 1 & 1022.47 & 2.456 \\
                & 1616.63 & 1 & 1000.94 & 2.512 & 1 & 1000.94 & 2.512 \\
  {\tt waterx}  & 1277.88 & 1 & 1024.85 & 10.848 & 7 & 997.27 & 71.696\\
  {\tt waterz}  & 1700.58 & - & - & - & - & - & - \\
                & 1497.95 & - & - & - & - & - & - \\
  \hline
\end{tabular}
\caption{Results obtained by applying the proposed local branching
  heuristic on instances taken from the MINLPLib. Time is expressed in
  seconds.}
\label{tab:results}
\end{table}

We tested the heuristic on 21 instances: for 18 instances we obtained
2 feasible solutions from {\tt couenne} within the 10 minutes time
limit, for the remaining 3 we only obtained 1 feasible solution. In
total, the heuristic was tested on 39 points used as incumbents. Note
that the instances: {\tt fo7}, {\tt fo7\_2}, {\tt fo8}, {\tt fo9},
{\tt m7}, {\tt o7}, {\tt o7\_2} are convex; therefore, the
linearization given by {\tt couenne} may not be the tightest, which
would be given by an outer approximation algorithm
\cite{fletcher_oa}. 

In 9 cases, the heuristic was not able to improve the incumbent within
the 10 iterations. On the {\tt waterz} instance, no feasible solution
was found by the algorithm. In 2 of the 9 unsuccessful runs (instances
{\tt fo7\_2, o7\_2}), it found a solution with the same objective
value as the incumbent; this may indicate the presence of symmetric
solutions. For 30 out of 39 initial feasible points ($76.9\%$),
Algorithm~\ref{alg:lb} was able to find a better solution. In 24 cases
($61.5\%$), an improved incumbent is found at the first iteration. The
improvement is significant. On the {\tt csched1}, {\tt csched1a}, {\tt
  ravem} and {\tt st\_e36} instances, our approach returns the best
known solution (reported on the MINLPLib website) from the first
feasible solution found by {\tt couenne}. On the {\tt o7\_2} instance,
the best known solution is returned from the second feasible solution
found by {\tt couenne}. The new incumbent is also very close to the
best known solutions for the instances: {\tt enprob48, enprob56,
  water4, waterx}. Large improvements are reported on the remaining
instances. Running time is typically less than 2.5 seconds; we remark
that we put a time limit of 2 seconds for Cplex. Therefore, the
running time can probably be reduced by decreasing Cplex's time
limit. On some instances ({\tt enpro48, o7, o7\_2, waterx}) several
seconds are required to solve the NLPs with {\tt ipopt}. This may be
due to lack of tuning of the parameters.

\section{Feasibility Heuristic}
An interesting observation is that our local branching heuristic found
at least one feasible solution on almost all the test instances. This
suggests employing a scheme similar to Algorithm~\ref{alg:lb} to
develop a primal heuristic whose purpose is only to find an initial
feasible point, regardless of its objective value. The first question
which arises is how to choose the point $x'$, which determines the
objective function for the MILP that is solved at the following
step. Since the purpose of this heuristic would be the discovery of a
first feasible solution, $x'$ should be a point in the interior of the
feasible region of ${\cal P}$. We can determine such point by solving
a continuous NLP over the feasible region of ${\cal P}$, with the
objective of maximizing the slacks between the constraints and their
respective bounds. For instance, if the problem is expressed in the
form (\ref{standard}), we should solve:
\begin{equation}
  \left. \begin{array}{rr\tS c\tS l} 
    \min_x \max_{j \in M} & g_j(x) && \\
    \forall j \in M & g_j(x) &\le& 0\\
    \forall i \in N & x_i^L \le x_i &\le& x_i^U \\
    \forall j \in N_I & x_j &\in& \mathbb{Z}, 
  \end{array} \right\} \mathcal{F}.
  \label{feas}
\end{equation}
This serves the purpose of finding a point $x'$ which is feasible and
maximizes the distance from the boundaries of the feasible region;
therefore, an integer feasible point near $x'$ is more likely to
satisfy the constraints. Since ${\cal F}$ may be a nonconvex problem,
depending on the constraints $g_j$'s, we solve it with a local NLP
solver in a multistart approach. Suppose we find $h$ local minima
$x'_1,\dots,x'_h$. Following the scheme of Algorithm~\ref{alg:lb}, the
next step is the solution of a MILP to obtain an integral feasible
solution to the convexification of the original problem, such that the
solution is near to one of the $x'_i, i=1,\dots,h$. This can be
modeled as a MILP. Clearly, since no feasible point is known, a local
branching constraint cannot be enforced when solving the
MILP. Therefore, the solution may take more time. However, we can use
early stopping criteria, such as a maximum time limit. Albeit this
approach is very simple, the computational experiments for the local
branching heuristic suggest that the idea might work.

\section{Conclusions and Future Research}
In this preliminary paper, we presented an idea for a local branching
heuristic that can be applied on nonconvex MINLPs with continuous and
binary variables. Our approach iteratively relies on solving a
sequence of MILPs and NLPs. We have reasons to believe that this
method is faster than to closely mimick the original idea of Fischetti
and Lodi \cite{localbranching} substituting the MILP solver with a
MINLP solver. Computational experiments run with a prototype of the
algorithm, written in AMPL, show that on most of the instances we are
able to significantly improve the incumbents, requiring a small CPU
time. We also observed that an approach similar to the one that we
presented could be used as initial feasibility heuristic, to be
employed at the beginning of the Branch-and-Bound tree. We did not
provide computational experiments for this idea. Our future research
will focus on the integration of the proposed techniques within an
existing solver for nonconvex MINLPs, to assess their usefulness in
practice.


\begin{thebibliography}{10}

\bibitem{couenne}
P.~Belotti.
\newblock Couenne, an open-source solver for mixed-integer nonconvex problems.
\newblock In preparation.

\bibitem{belottibbt}
P.~Belotti, J.~Lee, L.~Liberti, F.~Margot, and A.~W\"acther.
\newblock Branching and bounds tightening techniques for non-convex {MINLP}.
\newblock Technical Report RC24620, IBM, 2008.
\newblock \url{http://www.optimization-online.org/DB_HTML/2008/08/2059.html}

\bibitem{biegler}
L.~Biegler, I.~Grossmann, and A.~Westerberg.
\newblock {\em Systematic Methods of Chemical Process Design}.
\newblock Prentice Hall, Upper Saddle River (NJ), 1997.

\bibitem{bonmin}
P.~Bonami and J.~Lee.
\newblock {\tt BONMIN} user's manual.
\newblock Technical report, IBM Corporation, June 2007.

\bibitem{localbranching}
M.~Fischetti and A.~Lodi.
\newblock Local branching.
\newblock {\em Mathematical Programming}, 98:23--37, 2005.

\bibitem{fletcher_oa}
R.~Fletcher and S.~Leyffer.
\newblock Solving {M}ixed {I}nteger {N}onlinear {P}rograms by outer
  approximation.
\newblock {\em Mathematical Programming}, 66:327--349, 1994.

\bibitem{floudas}
C.~Floudas.
\newblock Global optimization in design and control of chemical process
  systems.
\newblock {\em Journal of Process Control}, 10:125--134, 2001.

\bibitem{vnsfoundations}
P.~Hansen and N.~Mladenovi\'{c}.
\newblock Variable neighbourhood search: Principles and applications.
\newblock {\em European Journal of Operations Research}, 130:449--467, 2001.

\bibitem{vnslb}
P.~Hansen, N.~Mladenovi\'c, and D.~Uro\v{s}evi\'c.
\newblock Variable neighbourhood search and local branching.
\newblock {\em Computers and Operations Research}, 33(10):3034--3045, 2006.

\bibitem{cplex110}
ILOG.
\newblock {\em ILOG CPLEX 11.0 User's Manual}.
\newblock ILOG S.A., Gentilly, France, 2007.

\bibitem{minlpbb}
S.~Leyffer.
\newblock User manual for {MINLP\_BB}.
\newblock Technical report, University of Dundee, UK, March 1999.

\bibitem{lln2}
L.~Liberti, C.~Lavor, and N.~Maculan.
\newblock A branch-and-prune algorithm for the molecular distance geometry
  problem.
\newblock {\em International Transaction in Operational Research}, 15:1--17,
  2008.

\bibitem{recipe}
L.~Liberti, G.~Nannicini, and N.~Mladenovi\'c.
\newblock A good recipe for solving {MINLP}s.
\newblock In V.~Maniezzo, T.~Stuetze, and S.~Voss, editors, {\em MATHEURISTICS:
  Hybridizing metaheuristics and mathematical programming}, Operations
  Research/Computer Science Interface Series. Springer, 2008.

\bibitem{nemhauser}
G.~Nemhauser and L.~Wolsey.
\newblock {\em Integer and Combinatorial Optimization}.
\newblock Wiley, New York, 1988.

\bibitem{tawarmalani1}
M.~Tawarmalani and N.~Sahinidis.
\newblock Global optimization of mixed integer nonlinear programs: A
  theoretical and computational study.
\newblock {\em Mathematical Programming}, 99:563--591, 2004.

\bibitem{ipopt}
A.~W\"achter and L.~T. Biegler.
\newblock On the implementation of a primal-dual interior point filter line
  search algorithm for large-scale nonlinear programming.
\newblock {\em Mathematical Programming}, 106(1):25--57, 2006.

\bibitem{wolsey}
L.~Wolsey.
\newblock {\em Integer Programming}.
\newblock Wiley, New York, 1998.

\end{thebibliography}

\end{document}